# Stochastic Dynamic Optimal Power Flow in Distribution Network with Distributed Renewable Energy and Battery Energy Storage

**Author information:**
Chenghui Tang, Jian Xu*, Yuanzhang Sun, Siyang Liao, Deping Ke, Xiong Li.
Department of Electrical Engineering, Wuhan University, Wuhan, Hubei, 430072, China
Correspondence and requests for materials should be addressed to C.K. (email: xujian@whu.edu.cn)

**Abstract:**

　　The penetration of distributed renewable energy (DRE) greatly raises the risk of distribution network operation such as peak shaving and voltage stability. Battery energy storage (BES) has been widely accepted as the most potential application to cope with the challenge of high penetration of DRE. To cope with the uncertainties and variability of DRE, a stochastic day-ahead dynamic optimal power flow (DOPF) and its algorithm are proposed. The overall economy is achieved by fully considering the DRE, BES, electricity purchasing and active power losses. The rainflow algorithm-based cycle counting method of BES is incorporated in the DOPF model to capture the cell degradation, greatly extending the expected BES lifetime and achieving a better economy. DRE scenarios are generated to consider the uncertainties and correlations based on the Copula theory. To solve the DOPF model, we propose a Lagrange relaxation-based algorithm, which has a significantly reduced complexity with respect to the existing techniques. For this reason, the proposed algorithm enables much more scenarios incorporated in the DOPF model and better captures the DRE uncertainties and correlations. Finally, numerical studies for the day-ahead DOPF in the IEEE 123-node test feeder are presented to demonstrate the merits of the proposed method. Results show that the actual BES life expectancy of the proposed model has increased to 4.89 times compared with the traditional ones. The problems caused by DRE are greatly alleviated by fully capturing the uncertainties and correlations with the proposed method.

# 0  Nomenclature

*Sets and indices*

| | | | |
|---|---|---|---|
| $i$ | Index of photovoltaic (PV) units. $i=1\ldots I$. | $n$ | Index of node. $n=1\ldots N$. |
| $j$ | Index of battery energy storage (BES). $j=1\ldots J$. | $s$ | Index of scenario of PV units actual power. $s=1\ldots S$ |
| $t$ | Index of time interval. $t=1\ldots T$. | | |

*Parameters*

| | | | |
|---|---|---|---|
| $s_{w,i}$ | Maximum apparent power capacity of the $i$-th PV unit | $SoC_j^{\max}$ | Upper bound of the SoC of $j$-th BES |
| $w_i^{\max}$ | Maximum active power capacity of the inverter of $i$-th PV unit | $SoC_j^{\min}$ | Lower bound of the SoC of $j$-th BES |
| $\alpha/\beta$ | Parameters of Beta distribution | $ch_j^{\max}$ | Upper limit on the rates of charge power of $j$-th BES |
| $\mu_i$ | Mean of Beta distribution of the $i$-th PV unit | $dis_j^{\max}$ | Upper limit on the rates of discharge of $j$-th BES |
| $\sigma_i^2 / \sigma_i$ | Variance/standard deviation of Beta distribution of the $i$-th PV unit | $\pi^s$ | Probability of $s$-th scenario |
| $i'$ | Imaginary part of the power | $V_0$ | Bus voltage magnitude of the substation |
| $c_{BES}$ | BES cell replacement price | $r_n / x_n$ | Resistance/reactance of branch $n$ |
| $Kch_j / Kdis_j$ | Number of charge/discharge half cycles of $j$-th BES. Indexed by $k$. | $\varepsilon$ | Deviation tolerance of the voltage magnitude |
| $k_1/k_2$ | Linear/nonlinear coefficients of the Polynomial cycle depth stress model | $c_{loss}$ | Cost coefficient of active power losses |
| $\underline{SoC}_{j,T}$ | Lower bound of end state of charge (SoC) of $j$-th BES | $c_{ss,p}$ | Cost coefficient of electricity purchasing |
| $\overline{SoC}_{j,T}$ | Upper bound of end SoC of $j$-th BES | $\rho$ | Permissible error in the multipliers update |
| $\eta_c / \eta_d$ | Charge/discharge efficiency of BES | $\gamma$ | Step size in the multipliers update |
| $\Delta t$ | Duration of each time interval | | |

*Variables*

| | | | |
|---|---|---|---|
| $w_{f,i}$ | Forecast power of $i$-th PV Units | $P_{0,t}$ | Active power from the transmission system |
| $w_{a,i}$ | Actual power of $i$-th PV Units | $Q_{0,t}$ | Reactive power from transmission system |
| $p_{j,t}^{ch}$ | Charge power of $j$-th BES at $t$ | | |
| $p_{j,t}^{dis}$ | Discharge power of $j$-th BES at $t$ | $d_{j,k}^{ch}$ | Cycle depths of $j$-th BES at $k$-th charge half cycle |
| $d_{j,k}^{ch}$ | Cycle depths of $j$-th BES at $k$-th charge half cycle | $q_{n+1,t}^{s(c)}$ | Consuming reactive power of node $n+1$ at $t$ of $s$-th scenario |
| $d_{j,k}^{dis}$ | Cycle depths of $j$-th BES at $k$-th discharge half cycle | $p_{n+1,t}^{s(g)}$ | Generating active power of node $n+1$ at $t$ of $s$-th scenario |
| $SoC_{j,t}$ | SoC of $j$-th BES at $t$. | $q_{n+1,t}^{s(g)}$ | Generating reactive power of node $n+1$ at $t$ of $s$-th scenario |
| $w_{a,i,t}^{s}$ | Actual active power of $s$-th scenario of $i$-th PV unit at $t$ | $V_{n,t}^{s}$ | Bus voltage magnitude of node $n$ at $t$ of $s$-th scenario |
| $q_{i,t}^{s}$ | Actual reactive power of $s$-th scenario of $i$-th PV unit at $t$ | $\lambda_{s,t,n}^{P}$ | Lagrange multiplier |
| $q_{i}^{s,\max}$ | Maximum reactive power capacity of $s$-th scenario of the $i$-th PV unit | $\lambda_{s,t,n}^{Q}$ | Lagrange multiplier |
| | | $\mu$ | Set of Lagrange multipliers |
| | | $d\mu$ | Set of Lagrange multipliers subgradient |
| $P_{n,t}^{s} / Q_{n,t}^{s}$ | Active/reactive power flows into the sending end of branch $n+1$ at $t$ under $s$-th scenario | | |
| $p_{n+1,t}^{s(c)}$ | Consuming active power of node $n+1$ at $t$ under $s$-th scenario | | |

*Functions*

| | | | |
|---|---|---|---|
| $F(\cdot)$ | Cumulative distribution function | | subproblem of distributed renewable energy (DRE) |
| $f(\cdot)$ | Probability density function | | |
| $f_{es,t}(\cdot)$ | Operation cost of BES at $t$ | $q^{BES}(\cdot)$ | Lagrange dual minimization problem of BES |
| $\Phi(d_{j,k}^{ch})$ | Cycle depth stress function | $q^{NP}(\cdot)$ | Lagrange dual minimization problem of node power |
| $\Phi(d_{j,k}^{dis})$ | Cycle depth stress function | $L(\cdot)$ | Lagrange dual function |
| $f_{loss,t}(\cdot)$ | Cost of active power losses at $t$ | $L^{W}(\cdot)$ | Subproblem of DRE function |
| $f_{ss,t}(\cdot)$ | Cost of electricity purchasing at $t$ | $L^{BES}(\cdot)$ | Subproblem of BES function |
| $q(\cdot)$ | Lagrange dual minimization problem | $L^{NP}(\cdot)$ | Subproblem of node power function |
| $q^{W}(\cdot)$ | Lagrange dual minimization | | |

*Abbreviations*

| | | | |
|---|---|---|---|
| DRE | Distributed renewable energy | CDF | Cumulative distribution function |
| BES | Battery energy storage | PDF | Probability density function |
| DOPF | Dynamic optimal power flow | SoC | State of charge |

# 1 Introduction

The high penetration of distributed renewable energy (DRE) involves new challenges around the efficiency and security of the distribution network. Battery energy storage (BES) has been widely accepted as the most potential application to cope with the challenge of DRE in the future power system. This paper focuses on the stochastic day-ahead dynamic optimal power flow (DOPF) in the distribution network with BESs and high penetration of DRE.

DOPF is an extension of optimal power flow to cover multiple time periods [1]. The day-ahead DOPF determines the charge/discharge power of BES, to cope with the challenges of peak shaving and voltage stability caused by DRE. Although there are many studies focus on it, two important aspects have not been adequately settled or even ignored. The first is the BES cell degradation cost, which is much larger than the maintenance cost and determines the real economy of the BES. The second is the uncertainties and correlations of DRE in the distribution network, which greatly affect the scheduled results of DOPF. Focusing on the above, this paper introduces three questions: how to model the BES and incorporate it into the DOPF; how to model the DRE and incorporate it into the DOPF; how to solve the day-ahead DOPF with BES and DRE. The literature review and the proposed method in this paper are discussed as follows.

**Literature review and the proposed methods**
**1) BES modeling**

A combined problem formulation for active-reactive optimal power flow in distribution networks with embedded wind generation and BES is proposed in [2]. However, the cost of BES is not considered in the above studies. Studies in [3]-[5] consider the investment cost in the siting and sizing problem of BES. Compared with the investment cost, cell degradation cost of BES is the real concern in the day-ahead DOPF. The degradation of BES is a complicated electrochemically process and related to many aspects [6]. The difficulty of the operation cost of BES is to build a model that is both accurate and easy for DOPF. Studies in [7] propose a physics-based degradation model for life prediction of Li-Ion batteries. Studies in [8] utilize electrochemical battery models to optimize the power management. However, it is difficult to involve the models in [7]-[8] into the DOPF. Studies in [9]-[11] model the cell degradation based on BES charge/discharge power. These studies would not capture the accumulation of previous charge and discharge effects of BES cell degradation and may severely deviate from the optimal operation.

Cycle-based cell degradation model can accurately model the fundamental degradation. The rainflow algorithm has been used for BES cell degradation model in [12]-[14]. However, it is difficult to incorporate the rainflow algorithm in the optimization problems because its procedure has no close form. To incorporate rainflow algorithm in the BES optimization, studies in [15] prove that the rainflow cycle-based cell degradation model is convex and provide a subgradient algorithm. However, studies in [15] optimize BES usage in frequency regulation market and not consider the real power system condition. In distribution network with BES, the objective function and constraints of BES are coupled with other components such as DER, load demand and transmission lines. Based on the convexity proof in [15], we further extend the rainflow cycle-based cell degradation model into the DOPF model. A stochastic day-ahead DOPF is proposed to achieve the overall economy by fully considering the DRE, BES, electricity purchasing and active power losses. The rainflow algorithm-based cycle counting method is incorporated in the proposed model to capture the cost of BES cell degradation.

**2) DRE modeling in the distribution network**

More and more DRE has been embedded in the distribution network. Studies in [1] incorporate energy storage and flexible demand in the active network management. A day-ahead economic dispatch of wind integrated power system is proposed considering demand response in [16]. However, renewable energy uncertainties have not been considered in the above studies [1][16]. To consider the renewable power uncertainties, Beta distribution is widely used in many studies [17]-[19]. Different from the concentrated wind power plants, DRE are involved in different nodes in the distribution network. Even if we can obtain their distribution models, it is difficult to directly employ them in the DOPF model.

Scenario-based method has been proved to be a good approach to cope with optimization problems in power system with DRE. Renewable power uncertainties and correlations can be represented accurately by renewable power scenarios. Studies in [20] generate photovoltaic (PV) units power scenarios by assuming that each PV unit follows a Beta distribution. However, PV units power in distribution network has strong dependence and this could not be considered by the method in [20]. Up to now, there are a great deal of studies [22]-[24] that focus on the wind power scenario generation methods in the transmission network while less of them are employed in the distribution network. Among all, Copula theory [24] shows a novel representation of the uncertainties and correlations for multiple wind power plants. To this end, we employ Copula theory for the scenario generation to capture the uncertainties and correlations of DRE in the distribution network.

**3) Solving the day-ahead DOPF with DRE scenarios**

Voltage instability caused by DRE brings an urgent need for the understanding of the uncertainties and correlations. However, although sufficient quantity of DRE power scenarios can be generated by the above methods. For the computation limit of the various problems such as DOPF, unit commitment and economic dispatch, large number of renewable power scenarios could not be employed in the optimization problem. Renewable power scenarios must be reduced to a smaller set. For instance, in studies [20], 1000 PV unit power scenarios are generated and reduced to 7 representative scenarios for stochastic optimal power flow using the fast forward reduction method. In studies [24], 3000 wind power scenarios are generated and reduced to a computationally feasible number 20 for unit commitment using the K-means clustering technique. In studies [25], 1000 wind power scenarios are generated and reduced to 20 scenarios for distribution feeder reconfiguration problem.

All the above scenario reduction methods aim to have a good approximation to the initial scenarios by the rest scenarios. However, with the increase of DRE embedding in the distribution network, the number of scenarios that incorporated in the optimization problem would be much limited. The uncertainties and correlations of DRE could not be fully considered by the traditional scenario reduction methods. To this end, we propose a Lagrange relaxation-based algorithm to incorporate much more scenarios in the DOPF model. The proposed algorithm has a significantly reduced complexity with respect to the existing scenario-based techniques and the size of the DOPF model increases with the number of scenarios by linear level.

In summary, to cope with the uncertainties and variability of DRE, we propose a stochastic day-ahead DOPF in distribution network. BES copes with the challenge of the peak shaving and voltage stability caused by high penetration of DRE.

In light of the above, the contributions of this paper are summarized as follows:

1) We propose a stochastic day-ahead dynamic optimal power flow of the distributed network to achieve the overall economy by fully considering the distributed renewable energy, battery energy storage, electricity purchasing and active power losses, in which the uncertainties and correlations of distributed renewable energy are fully considered based on scenarios.

2) The rainflow algorithm-based cycle counting method of battery energy storage is involved in the dynamic optimal power flow model to capture the cell degradation, which greatly extends the expected battery energy storage lifetime and achieve better economy.

3) We propose a Lagrange relaxation-based algorithm to solve the stochastic day-ahead dynamic optimal power flow. The proposed algorithm enables much more distributed renewable energy power scenarios incorporated in the optimal power flow model and better captures the uncertainties and correlations compared with the traditional ones.

The organization of the rest of this paper is as follows: Section 2 introduces the scenario generation method to consider the uncertainties and correlations of DRE. Section 3 describes the stochastic day-ahead DOPF model of the distributed network with DRE and BES. Section 4 proposes a Lagrange relaxation-based algorithm to solve the DOPF problem with a large number of DRE scenarios. Section 5 uses a case study for the IEEE 123-node test feeder to demonstrate the superiority and effectiveness of the proposed method. Finally, conclusions are drawn and future works are presented in Section 6.

# 2 Scenario Generation of DRE

## 2.1 Classical Model of DRE Actual Power

Conditional distribution models have been proved to be more accurately in considering the renewable energy uncertainties. Based on the conditional distributions, distributed renewable energy (DRE) power scenarios could be generated to show the possible occurrence of actual power. In this paper, we mainly discuss the photovoltaic (PV) units since they are more widely installed in the distribution network. Note that other types of DRE power plants such as wind power plants could also use the proposed model.

Beta distribution is a widely used model to represent the uncertainties of PV Units [17]-[20]. Assuming the maximum apparent power capacity of the $i$-th PV unit is $s_{w,i}$. The relation between $s_{w,i}$ and the maximum active power capacity of the inverter $w_i^{\max}$ is $w_i^{\max} = s_{w,i}/1.1$ [20]. Assuming the forecast power of PV Units is $w_{f,i}$, the occurrence of the actual power of PV Units $w_{a,i}$ is modeled by a Beta distribution as follows:

$$f_{Beta}(w_{a,i}) = \frac{1}{w_i^{\max}} (\frac{w_{a,i}}{w_i^{\max}})^{\alpha-1} \cdot (1 - \frac{w_{a,i}}{w_i^{\max}})^{\beta-1} \tag{1}$$

The following equations determine the relations between $\alpha$ and $\beta$ with the mean value $\mu$

and the variance $\sigma^2$ of the Beta distribution function:

$$\mu_i = w_{f,i} = \frac{\alpha}{\alpha+\beta} w_i^{max} \quad (2)$$

$$\sigma_i^2 = \frac{\alpha\beta}{(\alpha+\beta)^2 \cdot (\alpha+\beta+1)} (w_i^{max})^2 \quad (3)$$

where $\mu_i$ and $\sigma_i^2$ are the mean and variance parameters, respectively.

The relation between the forecast power of PV Units $w_{f,i}$ and the corresponding standard deviation of the Beta distribution function is as follows [17]-[20]:

$$\frac{\sigma_i}{w_i^{max}} = 0.2 \frac{w_{f,i}}{w_i^{max}} + 0.21 \quad (4)$$

The uncertainties of actual power of PV Units can be considered and PV Units power scenarios can be generated by the above model. However, the correlations between different PV Units could not be considered by the models in (1)~(4) since the power scenarios of each PV Unit are generated separately. The actual power of PV Units in the distributed network has strong correlation between each other due to the close geographical locations. Fig. 1 shows the joint distribution of a typical group of two PV Units power scenarios based on the models in (1)~(4) while Fig. 2 shows the joint distribution of their statistical historical PV Units power occurrence of the same forecast power [21]. We can see that although the uncertainties of each PV Unit can be considered by the above model, the PV Units power scenarios are far different with the occurrence of actual PV Units power. In order to consider both the uncertainties and correlations, the joint distribution of actual power of PV Units is required.

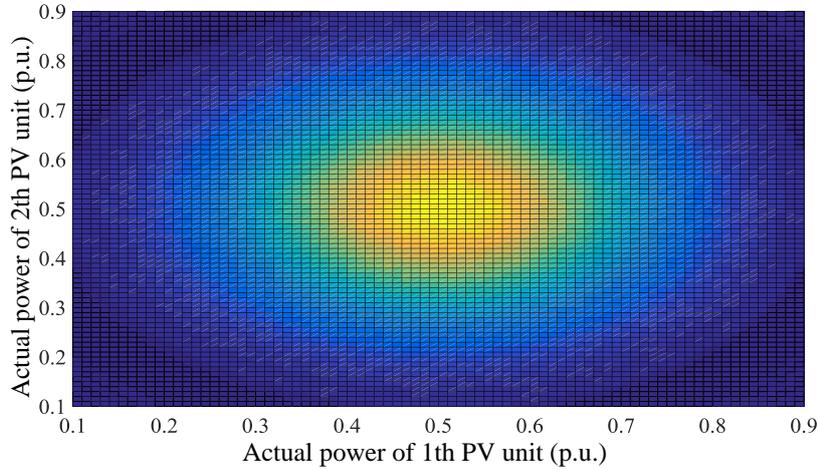

Figure 1. Joint distribution of a typical group of two PV Units power scenarios based on Beta distribution.

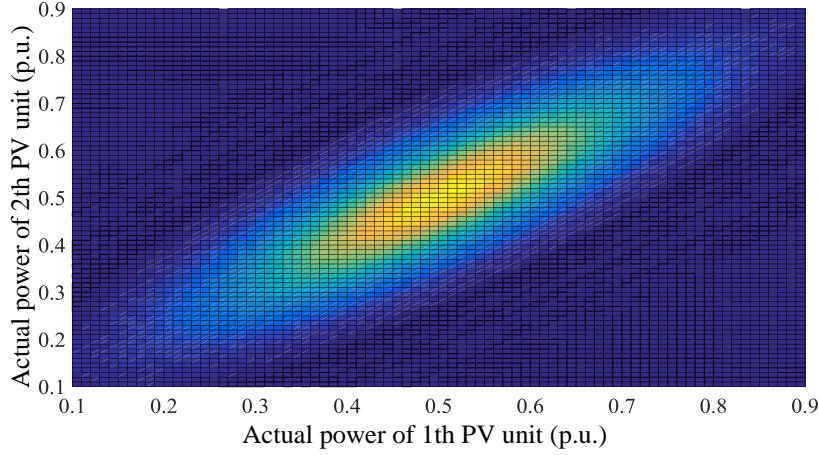

Figure 2. Joint distribution of statistical historical occurrence of actual PV Units power.

## 2.2 Conditional Joint Distribution Model of DRE Actual Power

Copula theory [24] is a good approach to consider the correlations between different PV units. Based on Copula theory, the joint cumulative distribution function (CDF) of actual power of PV Units $w_{a,i}$ and the forecast power of PV Units $w_{f,i}$ can be modeled as follows:

$$\begin{aligned}&F(w_{a,1}...w_{a,i}...w_{a,I}, w_{f,1}...w_{f,i}...w_{f,I})\\&=C(F(w_{a,1})...F(w_{a,i})...F(w_{a,I}), F(w_{f,1})...F(w_{f,i})...F(w_{f,I}))\end{aligned} \quad (5)$$

where $F(w_{a,i})$ is the marginal CDF of actual power of PV Units, $F(w_{f,i})$ is the marginal CDF of forecast power of PV Units.

In day-ahead dynamic optimal power flow (DOPF), the forecast power of each PV Unit can be obtained, then the conditional joint probability density function (PDF) is needed to model the actual power of all PV Units, as shown in (6).

$$\begin{aligned}&f(w_{a,1}...w_{a,i}...w_{a,I}|w_{f,1}...w_{f,i}...w_{f,I})\\&=\frac{c(F(w_{a,1})...F(w_{a,i})...F(w_{a,I}), F(w_{f,1})...F(w_{f,i})...F(w_{f,I}))}{c(F(w_{f,1})...F(w_{f,i})...F(w_{f,I}))} \cdot f(w_{a,1})...f(w_{a,i})...f(w_{a,I})\end{aligned} \quad (6)$$

Based on the joint conditional in (6), the correlations of the PV Units can be modeled. The uncertainties can be greatly reduced by considering the forecast power of each PV Units. To obtain the conditional joint PDF in (6), the first step is to find a suitable Copula function such as Gaussian copula, t-copula, empirical copula, etc and there are many studies on it. Which Copula form to choose is not the main concern in this paper, we use the most widely used Gaussian Copula and the method can be easy employed in other Copula form. Detailed information about Copula function could be found in [24]. Based on the Copula function in (6), we can generate PV units power scenarios that consider the uncertainties and correlations accurately.

# 3 DOPF Model in Distribution Network with High Penetration of DRE

## 3.1 Distribution Network Configuration

Fig. 3 shows a radical distribution network that includes substation, battery energy storages (BESs), DRE (PV units in this paper), and load demand. Assuming that there are $I$ PV units and $J$ BESs in the system. PV units offer the PV power to the system. The most of the power that offers to the system is from the substation. $P_0$ and $Q_0$ are the active power and reactive power that come from the transmission system, respectively.

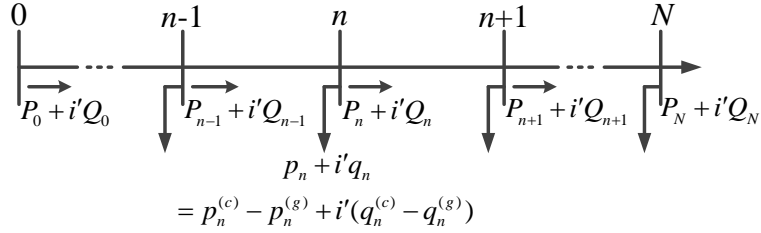

Figure 3. Diagram and nodes for the radial network.

## 3.2 BES Model

### 3.2.1 BES cell degradation model

BES can alleviate the problems brought about by the uncertainties and variability associated with DRE sources such as PV units [3][4]. Day-ahead DOPF determines the charge power $p_{j,t}^{ch}$ and discharge power $p_{j,t}^{dis}$ of BESs of each time interval in the next day. $t$ is the time interval (one hour in this paper), $t=1\ldots T$ and $T$ is the schedule horizon (24 hours). As we argued above, the cell degradation cost of BES is of vital importance to the day-ahead DOPF strategy. The BES cell degradation is captured by the rainflow cycle-based degradation model as follows:

$$f_{es}(d_{j,k}^{ch}, d_{j,k}^{dis}) = \sum_{j=1}^{J} f_{es,j}(d_{j,k}^{ch}, d_{j,k}^{dis}) = c_{BES} \sum_{j=1}^{J} (\sum_{k=1}^{Kch_j} \Phi(d_{j,k}^{ch}) + \sum_{k=1}^{Kdis_j} \Phi(d_{j,k}^{dis})) \qquad (7)$$

where $c_{BES}$ is the BES cell replacement price; $Kch_j$ and $Kdis_j$ are the number of charge and discharge half cycles of $j$-th BES, respectively; $d_{j,k}^{ch}$ and $d_{j,k}^{dis}$ are cycle depths of $j$-th BES and defined as below:

$$d_{j,1}^{ch}/d_{j,1}^{dis}\ldots d_{j,Kch}^{ch}/d_{j,Kdis}^{dis} = \text{Rainflow}(SoC_{j,1}\ldots SoC_{j,t}\ldots SoC_{j,T}) \qquad (8)$$

The rainflow cycle-based degradation model consists of (7) and (8) [15]. $\Phi(d_{j,k}^{ch})$ and $\Phi(d_{j,k}^{dis})$ are the cycle depth stress function, which defines the degradation of one half cycle under reference condition. $\sum_{k=1}^{Kch} \Phi(d_{j,k}^{ch})$ and $\sum_{k=1}^{Kdis} \Phi(d_{j,k}^{dis})$ calculate the total degradation by summing over all half cycles. $SoC_{j,t}$ is the state of charge (SoC) of $j$-th BES at time $t$.

### 3.2.2 Cycle depth stress function of BES

Cycle depth stress function $\Phi(d_{j,k}^{ch})$ and $\Phi(d_{j,k}^{dis})$ are critical part of the BES cell degradation model since it captures the cell aging caused by one (half) cycle under reference conditions [6]. There are different stress function forms for different types of batteries, such as linear cycle depth stress model [11], exponential cycle depth stress model [26] and Polynomial cycle depth stress model [27]. In this paper, the Polynomial cycle depth stress model is employed for it can capture the nonlinear impact of cycle depth on cell degradation under most conditions by lab tests. While other cycle depth stress models could also be employed in the proposed model. The Polynomial cycle depth stress model is as follows:

$$\Phi(d_{j,k}^{ch}) = k_1 \cdot (d_{j,k}^{ch})^{k_2}, \ \Phi(d_{j,k}^{dis}) = k_1 \cdot (d_{j,k}^{dis})^{k_2} \tag{9}$$

where $k_1$ and $k_2$ are model coefficients, which could be estimated by fitting battery cycling aging test data.

### 3.2.3 Others constraints of BES

Others constraints of BES are as follows:

$$SoC_{j,1} = SoC_{j,ini} \quad j \in J \tag{10}$$

$$\underline{SoC}_{j,T} \le SoC_{j,T} \le \overline{SoC}_{j,T} \quad j \in J \tag{11}$$

$$SoC_{j,t} = SoC_{j,t-1} + (p_{j,t}^{ch}\eta_c - p_{j,t}^{dis}/\eta_d) \cdot \Delta t \quad t \in T, j \in J \tag{12}$$

$$SoC_j^{\min} \le SoC_{j,t} \le SoC_j^{\max} \quad t \in T, j \in J \tag{13}$$

$$0 \le p_{j,t}^{ch} \le ch_j^{\max} \quad t \in T, j \in J \tag{14}$$

$$0 \le p_{j,t}^{dis} \le dis_j^{\max} \quad t \in T, j \in J \tag{15}$$

Equations (10) set the BES initial SoC in the day-ahead DOPF schedule horizon; Constraints (11) set the BES bounds of end SoC in the day-ahead DOPF schedule horizon; Equations (12) set the SoC at the end of time interval $t$ as a function of its SoC at the end of the previous time interval $t$-1 and the charge/discharge power that took place during time interval $t$; $p_{j,t}^{ch}$ and $p_{j,t}^{dis}$ are the charge and discharge power at time interval $t$ of $j$-th BES, respectively; $\eta_c$ and $\eta_d$ are the charge and discharge efficiencies of $j$-th BES, respectively; $\Delta t$ is the duration of each time interval. Constraints (13) impose an upper bound $SoC_j^{\max}$ and lower bound $SoC_j^{\min}$ on the SoC of $j$-th BES. Constraints (14) and (15) impose an upper bound $ch_j^{\max}$ and $dis_j^{\max}$ on the rates of charge and discharge power of $j$-th BES, respectively.

## 3.3 PV Units Power Model

PV units are connected on some user nodes in the distributed network. In this paper, actual power of the PV unit is modelled by the proposed scenario generation method in Section II. It is assumed that the scenario of $i$-th PV unit actual power at time $t$ is denoted by $w_{a,i,t}^s$. $s=1\ldots S$ and $S$ is the number of scenarios of PV units actual power. One scenario consists of the actual power scenarios of all PV units at same time interval. Each scenario has a probability $\pi^s$.

PV inverters translate the direct-current power to alternating current power. These PV inverters are also capable of generating or consuming reactive power by themselves [28]. It is assumed that the reactive power of *i*-th PV unit at time *t* of *s*-th scenario is denoted by $q_{i,t}^s$:

$$-q_i^{s,\max} \leq q_{i,t}^s \leq q_i^{s,\max} \quad i \in I, t \in T, s \in S \tag{16}$$

where $q_i^{s,\max} = \sqrt{s_{w,i}^2 - (w_{a,i,t}^s)^2}$ and $s_{w,i}$ is the maximum apparent power capacity of the *i*-th PV unit.

## 3.4 Power Flow Model

The classical linear power flow equations known as *LinDistFlow* [29]-[31] are employed as follows:

$$P_{n+1,t}^s = P_{n,t}^s - p_{n+1,t}^{s(c)} + p_{n+1,t}^{s(g)} \quad n \in N, t \in T, s \in S \tag{17}$$

$$Q_{n+1,t}^s = Q_{n,t}^s - q_{n+1,t}^{s(c)} + q_{n+1,t}^{s(g)} \quad n \in N, t \in T, s \in S \tag{18}$$

$$V_{n+1,t}^s = V_{n,t}^s - (r_n P_{n,t}^s + x_n Q_{n,t}^s)/V_0 \quad n \in N, t \in T, s \in S \tag{19}$$

$$(1-\varepsilon)^2 \leq \frac{V_{n,t}^s}{V_0} \leq (1+\varepsilon)^2 \quad n \in N, t \in T, s \in S \tag{20}$$

where $P_{n,t}^s$ and $Q_{n,t}^s$ are the active and reactive power flows into the sending end of branch *n*+1 connecting node *n* and node *n*+1 at time *t* under *s*-th DRE power scenario, respectively, as shown in Fig. 3; $p_{n+1,t}^{s(c)}$ and $q_{n+1,t}^{s(c)}$ are the consuming active and reactive power of node *n*+1 at time *t* under *s*-th DRE power scenario, respectively; $p_{n+1,t}^{s(g)}$ and $q_{n+1,t}^{s(g)}$ are the generating active and reactive power of node *n*+1 at time *t* under *s*-th DRE power scenario, respectively; $V_{n,t}^s$ is the bus voltage magnitude of node *n* at time *t* under *s*-th DRE power scenario; $V_0$ is the bus voltage magnitude of the substation and is assumed to be constant; $r_n$ and $x_n$ are the resistance and reactance of branch *n*; $\varepsilon$ is the deviation tolerance of the voltage magnitude of node *n* and is an user-defined value.

By the linearized model, the power losses in the distributed network is as follows:

$$P_{loss,t} = \sum_{s=1}^{S}[\pi^s P_{loss,t}^s] = \sum_{s=1}^{S}[\pi^s \sum_{n=0}^{N-1} r_n \frac{(P_{n,t}^s)^2 + (Q_{n,t}^s)^2}{V_0^2}] \quad t \in T, s \in S \tag{21}$$

$$f_{loss,t} = c_{loss} P_{loss,t} \quad t \in T \tag{22}$$

where $c_{loss}$ is the cost coefficient of active power losses.

## 3.5 Objective Function and Constraints

The objective function of stochastic day-ahead DOPF is to determine the charge and discharge power of BES $p_{j,t}^{ch}$, $p_{j,t}^{dis}$ and minimize the overall cost of the distribution network, including the cost of BES depreciation, the cost of active power losses and cost of electricity purchasing from the transmission network as follows:

$$\min f = \sum_{t=1}^{T}(f_{es,t} + f_{loss,t} + f_{ss,t})$$
$$= c_{BES}\sum_{j=1}^{J}[\sum_{k=1}^{Kch}\Phi(d_{j,k}^{ch}) + \sum_{k=1}^{Kdis}\Phi(d_{j,k}^{dis})] + c_{loss}\sum_{t=1}^{T}\sum_{s=1}^{S}[\pi^s\sum_{n=0}^{N-1}(r_n\frac{(P_{n,t}^s)^2 + (Q_{n,t}^s)^2}{V_0^2})] \quad (23)$$
$$+ c_{ss,p}\sum_{t=1}^{T}\sum_{s=1}^{S}[\pi^s P_{0,t}^s]$$

The constraints are as follows:

$$\begin{cases} d_{j,1}^{ch}/d_{j,1}^{dis}...d_{j,Kch}^{ch}/d_{j,Kdis}^{dis} = \text{Rainflow}(SoC_{j,1}...SoC_{j,t}...SoC_{j,T}) \; j \in J \\ SoC_{j,1} = SoC_{j,ini} \; j \in J \\ \underline{SoC}_{j,T} \leq SoC_{j,T} \leq \overline{SoC}_{j,T} \; j \in J \\ SoC_{j,t} = SoC_{j,t-1} + (p_{j,t}^{ch}\eta_c - p_{j,t}^{dis}/\eta_d)\cdot\Delta t \; j \in J, t \in T \\ SoC_j^{\min} \leq SoC_{j,t} \leq SoC_j^{\max} \; j \in J, t \in T \\ 0 \leq p_{j,t}^{ch} \leq ch_j^{\max} \; j \in J, t \in T \\ 0 \leq p_{j,t}^{dis} \leq dis_j^{\max} \; j \in J, t \in T \\ -q_i^{s,\max} \leq q_{i,t}^s \leq q_i^{s,\max} \; i \in I, t \in T, s \in S \\ P_{n+1,t}^s = P_{n,t}^s - p_{n+1,t}^{s(c)} + w_{a,i,t}^s \; i \in I, n \in N \cap I, t \in T, s \in S \\ P_{n+1,t}^s = P_{n,t}^s - p_{n+1,t}^{s(c)} + p_{j,t}^{ch} - p_{j,t}^{dis} \; j \in J, n \in N \cap J, t \in T, s \in S \\ P_{n+1,t}^s = P_{n,t}^s - p_{n+1,t}^{s(c)} \; n \in N \,|\, I,J, t \in T, s \in S \\ Q_{n+1,t}^s = Q_{n,t}^s - q_{n+1,t}^{s(c)} + q_{i,t}^s \; i \in I, n \in N \cap I, t \in T, s \in S \\ Q_{n+1,t}^s = Q_{n,t}^s - q_{n+1,t}^{s(c)} \; n \in N \,|\, I, t \in T, s \in S \\ V_{n+1,t}^s = V_{n,t}^s - (r_n P_{n,t}^s + x_n Q_{n,t}^s)/V_0 \; n \in N, t \in T, s \in S \\ (1-\varepsilon)^2 V_0 \leq V_{n,t}^s \leq (1+\varepsilon)^2 V_0 \; n \in N, t \in T, s \in S \end{cases} \quad (24)$$

where $N \cap I$ is referred to as the set of nodes that connected on the PV units; $N \cap J$ is referred to as the set of nodes that connected on the BES; $N \,|\, I,J$ is referred to as the set of nodes that do not connected on the PV units or BES; $N \,|\, I$ is referred to as the set of nodes that do not connected on the PV units. To clearly show the DOPF model, we assume that the PV unit and BES connect on different node of the distribution network. While the situation that PV unit and BES connect on the same node could also be considered by sample revision on the proposed model.

We can see from the model (23)~(24) that day-ahead DOPF is a very complicate optimization problem. The solving difficulty lies in the following three aspects：

- The model of BES has been a complicate problem while the cost and constraints of BES are coupled with other elements of the distribution network (DRE, nodes in the distribution network). For instance, we can use different models of BES (such as rainflow cycle-based degradation model, no operating cost model and linear power-based model). This would greatly influent the model type in (23)~(24).
- The model of DRE, active power and reactive power of each node, time interval are

- coupled with each other. This makes the size of the problem much larger.
- ● The size of the day-ahead DOPF increases with the number of DRE scenarios by index level. For this reason, the DRE scenarios need to be reduced to a very limited number, which could not accurately represent the renewable power uncertainties and correlations.

# 4 Solving the Day-Ahead DOPF by Lagrange Relaxation-Based Algorithm

In this section, we propose a Lagrange relaxation-based algorithm to solve the day-ahead DOPF with a large number of DRE power scenarios, which has a significantly reduced complexity with respect to the existing scenario-based techniques.

## 4.1 Lagrange Relaxation Problem

When we relax the transmission constraints (17) and (18), the Lagrange function of (23) ~(24) is as follows:

$$L(q_{i,t}^s, p_{j,t}^{ch}, p_{j,t}^{dis}, P_n^s, Q_n^s, \lambda_{s,t,n}^P, \lambda_{s,t,n}^Q)$$
$$= c_{BES} \sum_{j=1}^{J} [\sum_{k=1}^{Kch} \Phi(d_{j,k}^{ch}) + \sum_{k=1}^{Kdis} \Phi(d_{j,k}^{dis})] + c_{loss} \sum_{t=1}^{T} \sum_{s=1}^{S} [\pi^s \sum_{n=0}^{N-1} (r_n \frac{(P_{n,t}^s)^2 + (Q_{n,t}^s)^2}{V_0^2})]$$
$$+ c_{ss,p} \sum_{t=1}^{T} \sum_{s=1}^{S} [\pi^s P_{0,t}^s] - \sum_{s=1}^{S} \sum_{t=1}^{T} \sum_{n=0}^{N} \lambda_{s,t,n}^P (P_{n+1,t}^s - P_{n,t}^s + p_{n+1,t}^{s(c)}) + \sum_{s=1}^{S} \sum_{t=1}^{T} \sum_{i=1}^{I} \lambda_{s,t,n}^P w_{a,i,t}^s \quad (25)$$
$$+ \sum_{s=1}^{S} \sum_{t=1}^{T} \sum_{j=1}^{J} \lambda_{s,t,n}^P (p_{j,t}^{ch} - p_{j,t}^{dis}) - \sum_{s=1}^{S} \sum_{t=1}^{T} \sum_{n=0}^{N} \lambda_{s,t,n}^Q (Q_{n+1,t}^s - Q_{n,t}^s + q_{n+1,t}^{s(c)}) + \sum_{s=1}^{S} \sum_{t=1}^{T} \sum_{i=1}^{I} \lambda_{s,t,n}^Q q_{i,t}^s$$

Subject to:

$$\begin{cases} d_{j,1}^{ch}/d_{j,1}^{dis}...d_{j,Kch}^{ch}/d_{j,Kdis}^{dis} = \text{Rainflow}(SoC_{j,1}...SoC_{j,t}...SoC_{j,T}) \quad j \in J \\ SoC_{j,1} = SoC_{j,ini} \quad j \in J \\ \underline{SoC}_{j,T} \leq SoC_{j,T} \leq \overline{SoC}_{j,T} \quad j \in J \\ SoC_{j,t} = SoC_{j,t-1} + (p_{j,t}^{ch}\eta_c - p_{j,t}^{dis}/\eta_d) \cdot \Delta t \quad j \in J, t \in T \\ SoC_j^{\min} \leq SoC_{j,t} \leq SoC_j^{\max} \quad j \in J, t \in T \\ 0 \leq p_{j,t}^{ch} \leq ch_j^{\max} \quad j \in J, t \in T \\ 0 \leq p_{j,t}^{dis} \leq dis_j^{\max} \quad j \in J, t \in T \\ -q_i^{s,\max} \leq q_{i,t}^s \leq q_i^{s,\max} \quad i \in I, t \in T, s \in S \\ V_{n+1,t}^s = V_{n,t}^s - (r_n P_{n,t}^s + x_n Q_{n,t}^s)/V_0 \quad n \in N, t \in T, s \in S \\ (1-\varepsilon)^2 V_0 \leq V_{n,t}^s \leq (1+\varepsilon)^2 V_0 \quad n \in N, t \in T, s \in S \end{cases} \quad (26)$$

where $\lambda_{s,t,n}^P$ and $\lambda_{s,t,n}^Q$ is unbounded variables. The Lagrangian dual problem is to maximize the dual function:

$$q(\lambda_{s,t,n}^P, \lambda_{s,t,n}^Q) = \min\{L\} \quad s.t. \ (26) \tag{27}$$

(25)~(26) can be converted to the following three models, which is represent by subproblem of DRE, subproblem of BES and subproblem of node power:

Subproblem of DRE is as follows:

$$\begin{aligned}&L^W(q_{i,t}^s, \lambda_{s,t,n}^P, \lambda_{s,t,n}^Q) \\ &= \sum_{s=1}^S \sum_{t=1}^T \sum_{i=1}^I [\lambda_{s,t,n}^P w_{a,i,t}^s + \lambda_{s,t,n}^Q q_{i,t}^s] \\ &s.t. \ -q_i^{s,\max} \leq q_{i,t}^s \leq q_i^{s,\max} \quad i \in I, t \in T, s \in S\end{aligned} \tag{28}$$

Subproblem of BES is as follows:

$$\begin{aligned}&L^{BES}(p_{j,t}^{ch}, p_{j,t}^{dis}, \lambda_{s,t,n}^P, \lambda_{s,t,n}^Q) \\ &= c_{BES} \sum_{j=1}^J [\sum_{k=1}^{Kch} \Phi(d_{j,k}^{ch}) + \sum_{k=1}^{Kdis} \Phi(d_{j,k}^{dis})] + \sum_{s=1}^S \sum_{t=1}^T \sum_{j=1}^J \lambda_{s,t,n}^P (p_{j,t}^{ch} - p_{j,t}^{dis}) \\ &= \sum_{j=1}^J [c_{BES} [\sum_{k=1}^{Kch} \Phi(d_{j,k}^{ch}) + \sum_{k=1}^{Kdis} \Phi(d_{j,k}^{dis})] + \sum_{t=1}^T \sum_{s=1}^S \lambda_{s,t,n}^P (p_{j,t}^{ch} - p_{j,t}^{dis})] \\ &s.t. \begin{cases} d_{j,1}^{ch}/d_{j,1}^{dis}...d_{j,Kch}^{ch}/d_{j,Kdis}^{dis} = \text{Rainflow}(SoC_{j,1}...SoC_{j,t}...SoC_{j,T}) \ j \in J \\ SoC_{j,1} = SoC_{j,ini} \ j \in J \\ \underline{SoC}_{j,T} \leq SoC_{j,T} \leq \overline{SoC}_{j,T} \ j \in J \\ SoC_{j,t} = SoC_{j,t-1} + (p_{j,t}^{ch} \eta_c - p_{j,t}^{dis}/\eta_d) \cdot \Delta t \ j \in J, t \in T \\ SoC_j^{\min} \leq SoC_{j,t} \leq SoC_j^{\max} \ j \in J, t \in T \\ 0 \leq p_{j,t}^{ch} \leq ch_j^{\max} \ j \in J, t \in T \\ 0 \leq p_{j,t}^{dis} \leq dis_j^{\max} \ j \in J, t \in T \end{cases}\end{aligned} \tag{29}$$

Subproblem of node power is as follows:

$$\begin{aligned}&L^{NP}(P_{n,t}^s, Q_{n,t}^s, \lambda_{s,t,n}^P, \lambda_{s,t,n}^Q) \\ &= \sum_{t=1}^T [c_{loss} \sum_{s=1}^S [\pi^s \sum_{n=0}^{N-1} r_n \frac{(P_{n,t}^s)^2 + (Q_{n,t}^s)^2}{V_0^2}] + c_{ss,p} \sum_{s=1}^S [\pi^s P_{0,t}^s]] \\ &\quad - \sum_{s=1}^S \sum_{t=1}^T \sum_{n=0}^N \lambda_{s,t,n}^P (P_{n+1,t}^s - P_{n,t}^s + p_{n+1,t}^{s(c)}) - \sum_{s=1}^S \sum_{t=1}^T \sum_{n=0}^N \lambda_{s,t,n}^Q (Q_{n+1,t}^s - Q_{n,t}^s + q_{n+1,t}^{s(c)}) \\ &= \sum_{s=1}^S \sum_{t=1}^T [c_{loss} \pi^s \sum_{n=0}^{N-1} r_n \frac{(P_{n,t}^s)^2 + (Q_{n,t}^s)^2}{V_0^2} + c_{ss,p} \pi^s P_{0,t}^s \\ &\quad - \sum_{n=0}^N \lambda_{s,t,n}^P (P_{n+1,t}^s - P_{n,t}^s + p_{n+1,t}^{s(c)}) - \sum_{n=0}^N \lambda_{s,t,n}^Q (Q_{n+1,t}^s - Q_{n,t}^s + q_{n+1,t}^{s(c)})] \\ &s.t. \begin{cases} V_{n+1,t}^s = V_{n,t}^s - (r_n P_{n,t}^s + x_n Q_{n,t}^s)/V_0 \ n \in N, t \in T, s \in S \\ (1-\varepsilon)^2 V_0 \leq V_{n,t}^s \leq (1+\varepsilon)^2 V_0 \ n \in N, t \in T, s \in S \end{cases}\end{aligned} \tag{30}$$

Then $q(\lambda_{s,t,n}^P, \lambda_{s,t,n}^Q)$ can be decomposed into $q^W(\lambda_{s,t,n}^P, \lambda_{s,t,n}^Q)$, $q^{BES}(\lambda_{s,t,n}^P, \lambda_{s,t,n}^Q)$ and $q^{BP}(\lambda_{s,t,n}^P, \lambda_{s,t,n}^Q)$. (27) can be converted to the following three models:

$$\begin{cases} q^W(\lambda_{s,t,n}^P, \lambda_{s,t,n}^Q) = \min\{L^W(q_{i,t}^s, \lambda_{s,t,n}^P, \lambda_{s,t,n}^Q)\} \\ q^{BES}(\lambda_{s,t,n}^P, \lambda_{s,t,n}^Q) = \min\{L^{BES}(p_{j,t}^{ch}, p_{j,t}^{dis}, \lambda_{s,t,n}^P, \lambda_{s,t,n}^Q)\} \\ q^{NP}(\lambda_{s,t,n}^P, \lambda_{s,t,n}^Q) = \min\{L^{NP}(P_n^s, Q_n^s, \lambda_{s,t,n}^P, \lambda_{s,t,n}^Q)\} \end{cases} \tag{31}$$

By relaxing the global coupled constraints, the original problem in (23)~(24) can be decomposed into independent subproblems $q^W(\lambda_{s,t,n}^P, \lambda_{s,t,n}^Q)$, $q^{BES}(\lambda_{s,t,n}^P, \lambda_{s,t,n}^Q)$ and $q^{NP}(\lambda_{s,t,n}^P, \lambda_{s,t,n}^Q)$. We can see that the dual function $q(\lambda_{s,t,n}^P, \lambda_{s,t,n}^Q)$ can be evaluated by summing subproblems $q^W(\lambda_{s,t,n}^P, \lambda_{s,t,n}^Q)$, $q^{BES}(\lambda_{s,t,n}^P, \lambda_{s,t,n}^Q)$ and $q^{NP}(\lambda_{s,t,n}^P, \lambda_{s,t,n}^Q)$. The process of solving the subproblems $q^W(\lambda_{s,t,n}^P, \lambda_{s,t,n}^Q)$, $q^{BES}(\lambda_{s,t,n}^P, \lambda_{s,t,n}^Q)$ and $q^{NP}(\lambda_{s,t,n}^P, \lambda_{s,t,n}^Q)$ is referred to as the subproblem stage. After $q(\lambda_{s,t,n}^P, \lambda_{s,t,n}^Q)$ is evaluated, the multipliers are updated, which is referred to as the master problem stage. Note that the rainflow algorithm has been proved to be a convex optimization [15], the proposed DOPF model (23)~(24) is also convex optimization. Based on the Slater's condition [32] and strong duality theory, the duality gap between the Lagrangian dual problem and the primary problem in (23)~(24) is zero. This means that by solving the subproblems and master problem, the proposed DOPF model (23)~(24) can be solved optimally.

## 4.2 Subproblem Solution

After updating the multipliers $\lambda_{s,t,n}^P, \lambda_{s,t,n}^Q$ (this will be further discussed in the subsection 4.3), the subproblems (31) are solved. We can see that (28)~(30) are independent from each other. Parallel computation can be employed in the sunproblems solution. The above three sunproblems are solved as follows.

### 4.2.1 Subproblem of DRE

We can see from (28) that the subproblem of DRE is simple linear programming problems. Each scenario *s*, each time interval *t* and each DRE are independent from each other. This means that the subproblem of DRE has very high computational efficiency and flexible allocation of computing ability by the parallel computation.

$$q^W(\lambda_{s,t,n}^P, \lambda_{s,t,n}^Q) = \sum_{s=1}^{S}\sum_{t=1}^{T}\sum_{i=1}^{I} q_{s,t,i}^W(\lambda_{s,t,n}^P, \lambda_{s,t,n}^Q)$$
$$= \sum_{s=1}^{S}\sum_{t=1}^{T}\sum_{i=1}^{I} \min\{L_{s,t,i}^W(q_{i,t}^s, \lambda_{s,t,n}^P, \lambda_{s,t,n}^Q)\}$$
(32)

### 4.2.2 Subproblem of BES

We can see from (29) that the subproblem of BES is independent from each BES. Parallel computation can be employed in the sunproblem solution of BES. What's more, from the item of $\sum_{s=1}^{S}\lambda_{s,t,n}^P(p_{j,t}^{ch} - p_{j,t}^{dis})$ in (29) we can see that the number of DRE scenarios only influent the coeefficience of $p_{j,t}^{ch} - p_{j,t}^{dis}$. This means that the size of subproblem of BES would not be increased with the number of DRE scenarios. This allows a high computational efficiency for each BES subproblem solution.

$$q^{BES}(\lambda_{s,t,n}^P, \lambda_{s,t,n}^Q) = \sum_{j=1}^{J} q_j^{BES}(\lambda_{s,t,n}^P, \lambda_{s,t,n}^Q)$$
$$= \sum_{j=1}^{J} \min\{L_j^{BES}(p_{j,t}^{ch}, p_{j,t}^{dis}, \lambda_{s,t,n}^P, \lambda_{s,t,n}^Q)\}$$
(33)

The single BES subproblem can be solved by the function named as "fmincon" in matlab.

For the reason that the model in (29) is convex [15], the global optimization result could be easy obtained.

### 4.2.3 Subproblem of node power in the distribution network

We can see from (30) that the subproblem of node power in the distribution network are independent from each scenario $s$ and each time interval $t$. Each subproblem of one scenario $s$ and one time interval $t$ are simple quadratic programming problem. Parallel computation can be employed in the subproblem solution of node power.

$$q^{NP}(\lambda_{s,t,n}^P, \lambda_{s,t,n}^Q) = \sum_{s=1}^{S}\sum_{t=1}^{T} q_{s,t}^{NP}(\lambda_{s,t,n}^P, \lambda_{s,t,n}^Q)$$
$$= \sum_{s=1}^{S}\sum_{t=1}^{T} \min\{L_{s,t}^{NP}(P_{n,t}^s, Q_{n,t}^s, \lambda_{s,t,n}^P, \lambda_{s,t,n}^Q)\} \quad (34)$$

## 4.3 Master Problem Solution

In this paper, we propose a subgradient method to update the multipliers. The set of $\lambda_{s,t,n}^P, \lambda_{s,t,n}^Q$ is defined by $\mu$, i.e., $\mu = \{\lambda_{s,t,n}^P, \lambda_{s,t,n}^Q \; \forall s,t,n\}$. The set of elements of the subgradient $d\lambda_{s,t,n}^P, d\lambda_{s,t,n}^Q$ is defined by $d\mu$, i.e., $d\mu = \{d\lambda_{s,t,n}^P, d\lambda_{s,t,n}^Q \; \forall s,t,n\}$, which are expressed as:

$$d\lambda_{s,t,n}^P = \begin{cases} P_{n,t}^s - p_{n+1,t}^{s(c)} + w_{a,i,t}^s - P_{n+1,t}^s & i \in I, n \in N \cap I, t \in T, s \in S \\ P_{n,t}^s - p_{n+1,t}^{s(c)} + p_{j,t}^{ch} - p_{j,t}^{dis} - P_{n+1,t}^s & j \in J, n \in N \cap J, t \in T, s \in S \\ P_{n,t}^s - p_{n+1,t}^{s(c)} - P_{n+1,t}^s & n \in N \mid I, J, t \in T, s \in S \end{cases} \quad (35)$$

$$d\lambda_{s,t,n}^Q = \begin{cases} Q_{n,t}^s - q_{n+1,t}^{s(c)} + q_{n+1,t}^{s(g)} + q_{i,t}^s - Q_{n+1,t}^s & i \in I, n \in N \cap I, t \in T, s \in S \\ Q_{n,t}^s - q_{n+1,t}^{s(c)} + q_{n+1,t}^{s(g)} - Q_{n+1,t}^s & n \in N \mid I, t \in T, s \in S \end{cases} \quad (36)$$

The procedure for solving the master problem is as follows.
1) Initialize the multipliers $\mu^{(k)} = \{\lambda_{s,t,n}^P{}^{(k)}, \lambda_{s,t,n}^Q{}^{(k)} \; \forall s,t,n\}$. Set the permissible error $\rho > 0$ and the step size $\gamma$. $\chi = 0$ in this step.
2) Solve the subproblems in (28)~(30) and calculate the subgradient $d\mu^{(k)}$ at $\mu^{(k)}$.
3) Check whether or not the convergence criterion is satisfied. If $\|d\mu^{(k)}\|_\infty < \rho$, terminate the algorithm; otherwise, update the multipliers by $\mu^{(k+1)} = \mu^{(k)} + \gamma \cdot d\mu^{(k)}$, set $\chi = \chi + 1$ and return to step 2).

## 4.4 Flowchart of the Proposed Lagrange Relaxation-Based Algorithm

Fig. 4 shows the flowchart of the proposed Lagrange relaxation-based algorithm. We can note that the computation time of the proposed algorithm is mainly determined by the subproblems. The subproblems of DRE, BES and node power have very flexible allocation of computing ability for the reason that each $\{q_{i,t}^s\}$, $\{p_{j,t}^{ch}, p_{j,t}^{dis}\}$ and $\{P_{n,t}^s, Q_{n,t}^s\}$ is independent from each other, as shown in the purple, blue and green boxes, respectively. Parallel computation could be employed and achieve very high computational efficiency. Another important aspect is that with the increase of number of DRE scenarios, the size of the day-ahead DOPF increases by linear level. Due to the above features, the proposed algorithm enables much more DRE power scenarios incorporated in the DOPF and better captures the uncertainties and correlations compared with solving the primary problem in (23)~(24).

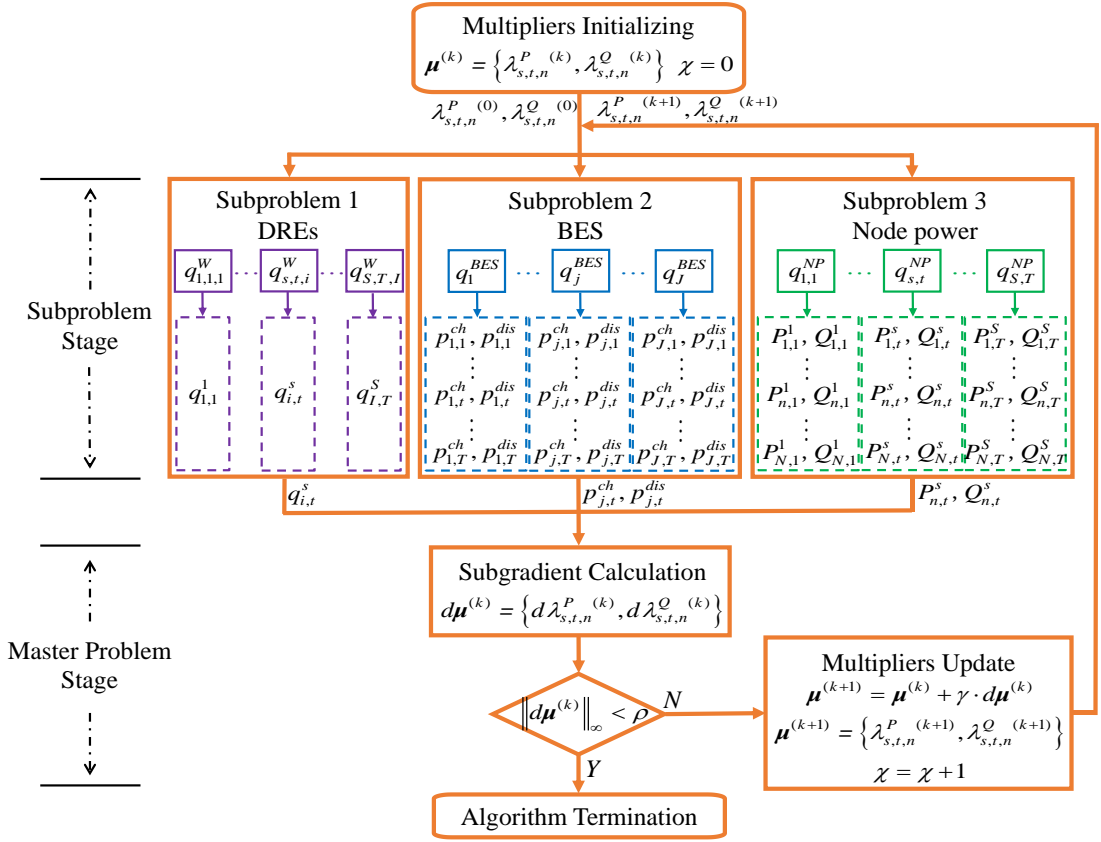

Figure 4. Flowchart of the proposed Lagrange relaxation-based algorithm.

# 5 Case Study

## 5.1 Basic Information and Data

Fig. 5 shows the topology of the IEEE 123-node test feeder [33]. The number of PV units is 8 and their capacities are given in Table. I. The data of PV units come from a real distribution network in Shandong Province, China. The penetration rate of DRE (PV units power) in this system is around 25.6% of the annual electricity consumption. 1000 PV units power scenarios are generated and incorporated in the proposed day-ahead DOPF model. The number of BES is 6 and each capacity is 75kWh, as shown in Table. II, as shown in Fig. 5. Polynomial cycle depth stress model is employed in this paper and $k_1$ and $k_2$ in (9) is $4.5\times 10^{-4}$ and 2.2, respectively. The allowance voltage offset $\varepsilon$ is 0.02. BES cell replacement price is 150$/kWh [34]. The charge and discharge efficiency $\eta_c, \eta_d$ are both 0.95. The BES has a 4h energy rating [19]. The initial SoC of all BES are both 0.6p.u.. The lower and upper bounds of end SoC of all BES are 0.3p.u. and 0.7p.u., respectively. The algorithm is run in MATLAB R2013a on a Core-i5 2.70-GHz notebook computer.

Table. I. Ratings of the DRE in the IEEE 123-node test feeder.

| Type | DRE-1 | DRE-2 | DRE-3 | DRE-4 | DRE-5 | DRE-6 | DRE-7 | DRE-8 |
|---|---|---|---|---|---|---|---|---|
| Node | 23 | 35 | 47 | 52 | 62 | 77 | 89 | 101 |
| Power Rating (kVA) | 450 | 450 | 450 | 450 | 450 | 300 | 450 | 300 |

Table. II. Ratings of the BES in the IEEE 123-node test feeder.

| Type | BES-1 | BES-2 | BES-3 | BES-4 | BES-5 | BES-6 |
|---|---|---|---|---|---|---|
| Node | 7 | 21 | 35 | 57 | 76 | 197 |
| Power Rating (kW) | 75 | 75 | 75 | 75 | 75 | 75 |

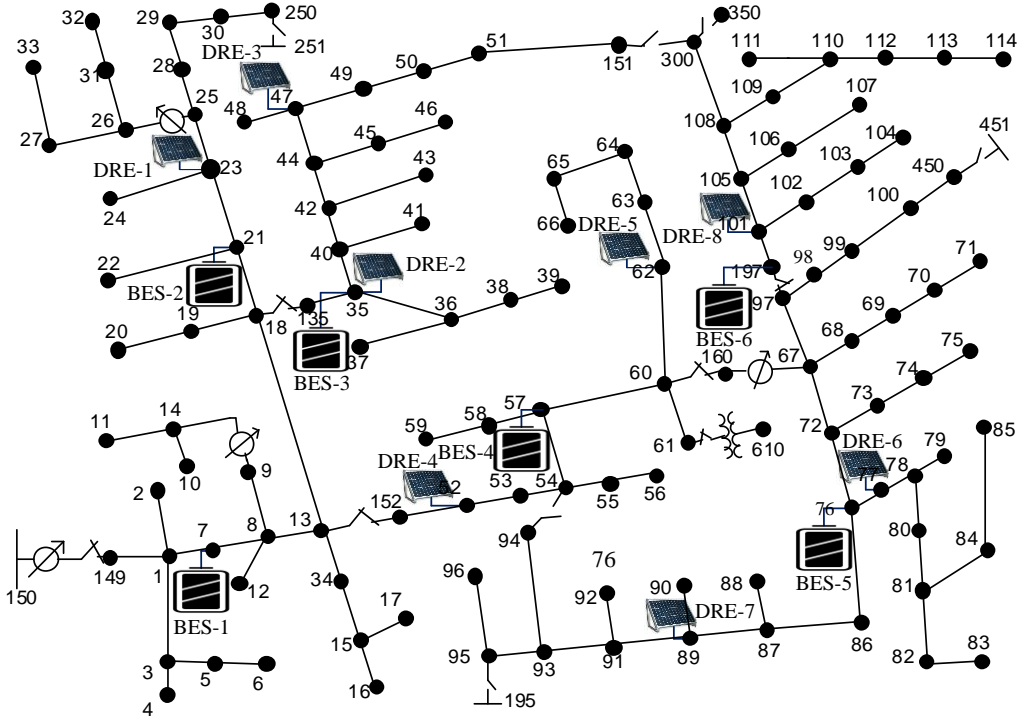

Figure 5. Topology of the IEEE 123-node test feeder with DRE and BES.

## 5.2 Charge/Discharge Power Based on Different BES Models

The proposed day-ahead DOPF model is solved based on three different BES cell degradation models: the proposed rainflow cycle-based cell degradation model, the no operating cost model [35] and linear power-based model [11]. For the reason that they are all convex, the above three BES degradation models can be easily incorporated into the proposed day-ahead DOPF model.

Fig. 6 shows the day-ahead electricity price and the SoC evolutions of BES-1 under the above three BES degradation models. We can see that the SoC of the three BES degradation models (brown, purple and red lines) have similar trends: charging at the time with lower electricity price (e.g., 00:00~02:00), discharging at the time with higher electricity price (e.g., 07:00~09:00, 19:00~21:00). The BES plays the role of peak shaving with the help of electricity price. An important thing need to be noticed is that the BES of three degradation models all charge between 11:00~14:00. The reason is that the PV units power are at the peak and raise the voltage of the distribution network. BES absorbs the PV units power by charging to deal with the potential voltage instability. Voltage instability caused by PV units power uncertainties would be further discussed in subsection 5.3 and 5.4.

Among all, SoC of no operating cost model (brown line) has deepest charge/discharge depths to gain the benefit of peak shaving without considering the cell degradation. Charge/discharge depths of linear power-based model are lower than them of no operating cost model. We use the algorithm in [36] to count the charge/discharge cycles. Compared with the linear power-based model (purple line), the proposed rainflow cycle-based cell degradation model (red line) avoids large charge/discharge circles, which could greatly accurate the BES cell degradation.

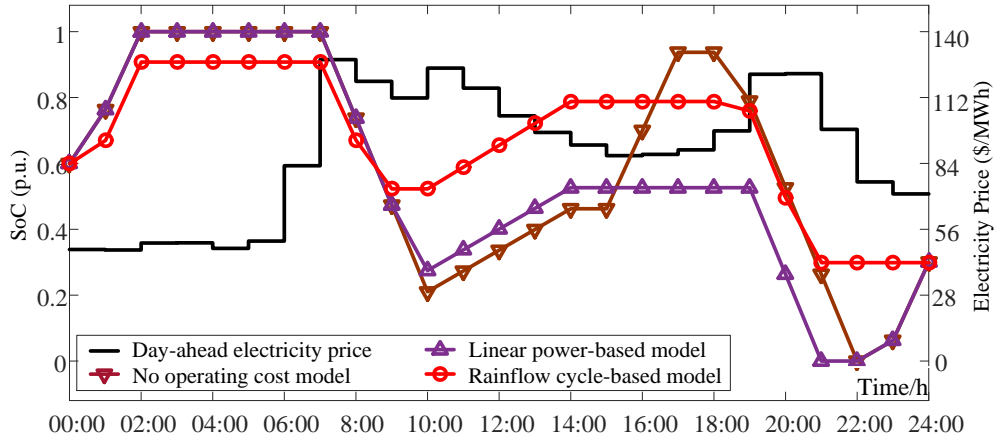

Figure 6. SoC evolutions of BES-1 and day-ahead electricity price curve.

Table. III shows the cost under the above three BES degradation models. We can see that the no operating cost model and linear power-based model all attempt to get the peak shaving benefit but failed to capture the cell degradation cost. As a result, actual net benefit of BES could not be guaranteed. We can see that the BES cell degradation cost is a relatively very high cost for the scheduling of BES and need to be fully considered. The proposed rainflow cycle-based degradation model greatly reduces the actual BES cell degradation cost. Compared with no operating cost model and linear power-based model, the actual BES life expectancy of the proposed rainflow cycle-based degradation model has increased to 4.89 times and 2.62 times, respectively.

Table. III. Costs of different BES degradation models.

| Average Daily Cost & BES life expectancy | No operating cost model | Linear power-based model | Rainflow cycle-based model |
|---|---|---|---|
| **Modeled BES cost ($)** | 0 | 4.93 | 15.06 |
| **Peak shaving benefit ($)** | 33.00 | 32.38 | 24.89 |
| **Modeled net benefit of BES ($)** | 33.00 | 27.45 | 9.83 |
| **Actual BES cost ($)** | 73.75 | 39.45 | 15.06 |
| **Actual net benefit of BES ($)** | -40.75 | -7.07 | 9.83 |
| **Actual BES life expectancy (year)** | 2.51 | 4.69 | 12.28 |

## 5.3 The Effect of the Uncertainties and Correlations of DRE

In this subsection, we analyze the effect of the uncertainties and correlations of DRE by comparing the day-ahead DOPF based on three DRE model, i.e., the deterministic DRE power model, i.e., does not consider DRE uncertainties, the scenario generation method without considering the DRE correlations [20] and the proposed DRE power scenario generation method. Based on the charge/discharge power obtained by the day-ahead DOPF models with the three DRE models, we generate another 5000 scenarios to test the system and get the voltage instability frequency (number of scenarios that cause voltage instability) of each DRE model, as shown in the Table IV. We can see that by considering the uncertainties and correlations of DRE, the voltage instability probability could be greatly reduced.

Table. IV. Voltage instability frequency of different DRE models.

|  | Proposed scenario generation model | Deterministic model | Scenario model without considering correlations |
|---|---|---|---|
| **Voltage instability frequency** | 4 | 421 | 163 |
| **Voltage instability probability** | 0.08% | 8.42% | 3.26% |

To clearly show the effect of the uncertainties and correlations, we get the sum power scenario of DRE by summing the active power of all PV units at same time interval, as shown in Fig. 7. What's more, we draw the boundary of certain numbers of sum power scenario to clearly show their distribution. For instance, 70% confidence level means the area of 70% scenarios located in, i.e., the deepest red area in Fig. 7(a) and the deepest green area in Fig. 7(b).

We can see from Fig. 7(a) and Fig. 7(b) that the variance of the distribution of sum power scenario and forecast power has nearly same tendency, i.e., increasing from morning, reach the peak at noon and then decrease to zero at night. This is in consistent with the Beta model in (4). Compared with the distributions in Fig. 7(a) and Fig. 7(b), smaller variance (less overall uncertainties) is detected if the DRE correlations are not considered. This means that overall uncertainties of DRE would be underestimated if we do not consider the DRE correlations. The reason is as follows.

The overall uncertainties of DRE are not simply related to the uncertainty of each PV unit, but also related to the correlations between all the PV units. We can take the joint Gaussian distribution of two PV units ($A$ and $B$) for instance. The variance of the sum PV units power $D(A+B) = D(A) + D(B) + 2Cov(A, B)$, where $Cov(A, B)$ is the covariance of $A$ and $B$.

Note that the PV units in the distribution network usually have strong correlations. The overall uncertainties could not be caught by separately generating scenarios of each PV unit as in [20]. The underestimation of overall uncertainties of DRE causes the underestimation of potential voltage instability. To this end, the BES have not offered enough charge power to maintain the voltage stability if the scenario model in [20] or deterministic DRE power model are employed.

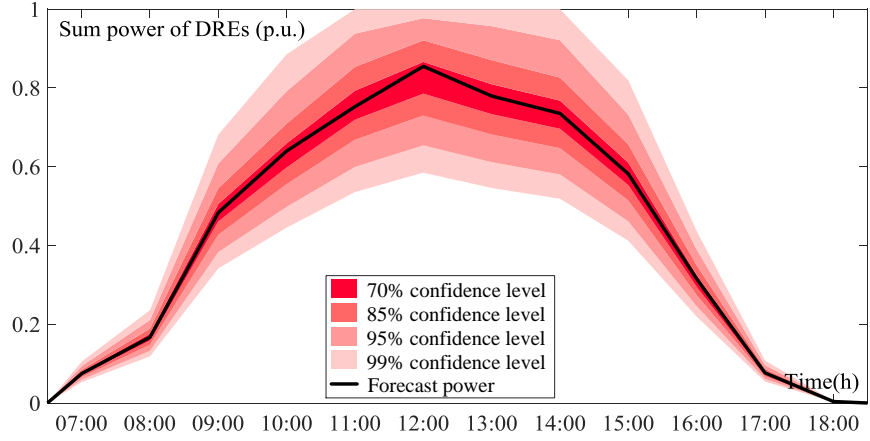

(a). Sum power scenarios of the proposed DRE power scenario generation method.

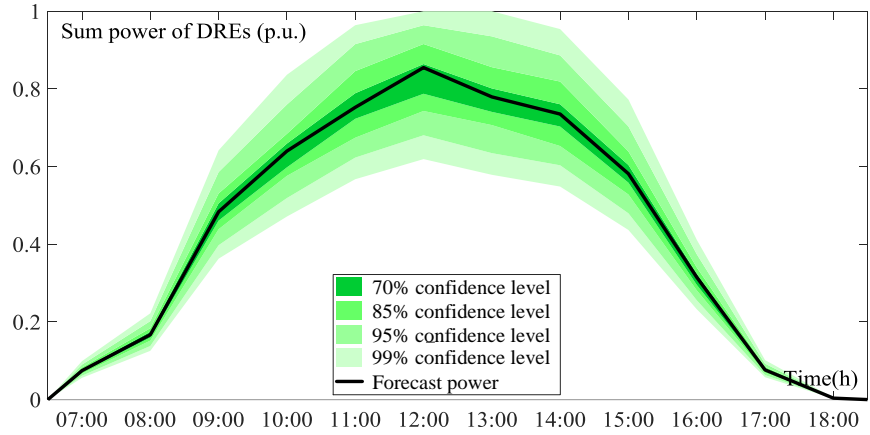

(b). Sum power scenarios of the method without considering the DRE correlations

Figure 7. Distribution of DRE scenarios based on different scenario generation method.

## 5.4 Comparison Between the Proposed Algorithm and Traditional Scenario Reduction Methods

In this subsection, we compare the proposed Lagrange relaxation-based algorithm and traditional scenario reduction-based method. Scenario reduction-based method solve the DOPF model in (23)~(24) with DRE scenarios obtained by the scenario reduction method in [37].

Similar to the subsection 5.3, we use 5000 scenarios to test the real system and get the voltage instability frequency of the scenario reduction-based model, as shown in the Table V. We can see that with the increase of the number of rest scenarios (7, 14 and 21), voltage instability could be improved. However, the voltage instability probability still maintains a high level compared with the proposed method. The reason is that traditional scenario reduction method try to have a good approximation to the initial scenarios by the rest scenarios. With the increase of DRE, the rest scenarios would be closer to the median area and underestimate the DRE uncertainties.

Table. V. Voltage instability frequency of scenario reduction-based method.

|  | Number of the rest scenario | | |
|---|---|---|---|
|  | 7 | 14 | 21 |
| **Voltage instability frequency** | 241 | 210 | 183 |
| **Voltage instability probability** | 4.82% | 4.20% | 3.66% |

As discussed above, the computation time of the proposed algorithm is primarily determined by the subproblems. Table VI shows the computation time of the subproblems of the proposed method. The computation time of subproblem of BES is relatively higher while it nearly not increase with the increase of scenario number. The subproblem of DRE and subproblem of node power all nearly increase with the increase of scenario number by linear level. Fig. 8 compared the computation time of the proposed algorithm and scenario reduction-based method with the same number of scenarios, i.e., 7, 14 and 21. We can see that the computation time of the proposed algorithm increases linear level with the scenario number and is much less time with the increase of scenarios. Scenario reduction-based method has competitive computational efficiency with 7 scenarios. However, computation times become even larger with the increase of scenario number. When the scenario number is 21, it is hard to solve in the day-ahead time period. To this end, the proposed algorithm enables much more scenarios incorporated in the DOPF model and better capture the DRE uncertainties and correlations.

Note that the computation times of the proposed algorithm in Fig. 8 are not the parallel computing. As argued above, the subproblems of DRE, BES and node power have very flexible allocation of computing ability for the reason that each $\{q_{i,t}^s\}$, $\{p_{j,t}^{ch}, p_{j,t}^{dis}\}$ and $\{P_{n,t}^s, Q_{n,t}^s\}$ is independent from each other. Parallel computation could be employed and achieve very high computational efficiency in real power system operation.

Table. VI. Computation times of the subproblem of the proposed method.

|  | Subproblem of DRE | Subproblem of BES | Subproblem of node power |
|---|---|---|---|
| **Time (s)** | ≤0.001 | 3.16 | 0.14 |

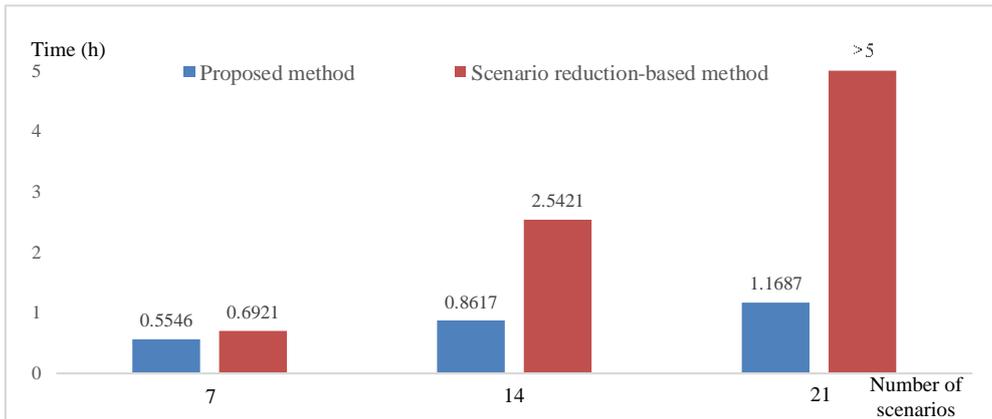

Figure 8. Computation times of the proposed method and scenario reduction-based method.

# 6  Conclusions and Future Works

This paper proposes a stochastic DOPF in distribution network with DRE and BES. The rainflow cycle-based degradation model is incorporated in the proposed model. DRE uncertainties and correlations are captured by the Copula theory and fully considered by the proposed Lagrange relaxation-based algorithm. Results show that by involving the rainflow cycle-based degradation model, the actual BES life expectancy has increased to 4.89 times compared with traditional no operating cost model. Voltage instability probability caused by high DRE penetration is greatly reduced by the BES scheduling and full consideration of the uncertainties and correlations of DER. The proposed Lagrange relaxation-based algorithm has a significantly reduced complexity with respect to traditional scenario reduction-based method and better capture the DRE uncertainties and correlations.

BES has been widely accepted as the most potential application to cope with the challenge of high penetration of DRE. However, the flexible application of BES is constrained by the present price. With the development of the BES technique especially the reduction of BES cell price and the increase of full charge/discharge cycle number, BES would play an important part of the power system for various usages.

Future works on DOPF in distribution network with DRE and BES would be focused on the following three aspects: BES operation in distribution network with real-time regulation market; considering the correlation of DER and the load demand; modeling the demand response and incorporating it into the DOPF model. BES could gain more benefit under the real-time regulation market under the present cell price. The number of load demand in the distribution network is much larger than the number of DRE, which would bring dimension disaster if the scenario method is employed. In this premise, how to incorporate demand response into the DOPF model is a further problem. The above problems would be further studied in future works.

**Acknowledgements**

This work was supported in part by the National Key Research and Development Program of China (2016YFB0900105), in part by the National Natural Science Foundation of China (51477122)